\documentclass[12pt,a4paper]{article}
\usepackage{geometry,amssymb,latexsym,amsmath,graphicx,amsthm}
\theoremstyle{plain}
\def \ni{\noindent}

\newtheorem{theorem}{Theorem}[section]
\newtheorem{lem}[theorem]{Lemma}
\newtheorem{definition}[theorem]{Definition}
 
\newtheorem{coro}{Corollary} 
\newtheorem{case}{Case} \newtheorem{rem}{Remark}[section] \newtheorem{cla}{Claim}

\newtheorem{obs}[theorem]{Observation}
\newtheorem{prop}[theorem]{Proposition}

\newcommand{\bl}{\begin{lem}}
\newcommand{\el}{\end{lem}}
\newcommand{\bt}{\begin{theorem}}
\newcommand{\et}{\end{theorem}}
\newcommand{\bc}{\begin{coro}}
\newcommand{\ec}{\end{coro}}
\newcommand{\bd}{\begin{definition}}
\newcommand{\ed}{\end{definition}}
\newcommand{\bp}{\begin{proof}}
\newcommand{\ep}{\end{proof}}
\newcommand{\bo}{\begin{obs}}
\newcommand{\eo}{\end{obs}}
\newcommand{\bpr}{\begin{prop}}
\newcommand{\epr}{\end{prop}}
\newcommand{\br}{\begin{rem}}
\newcommand{\er}{\end{rem}}
\newcommand{\bcl}{\begin{cla}}
\newcommand{\ecl}{\end{cla}}
\newcommand{\bca}{\begin{case}}
\newcommand{\eca}{\end{case}}
\def \ni{\noindent}

\newcommand{\gpk}{\gamma_{{\rm P},k}}		

\newcommand{\kPDS}{{$k$-PDS}}

\newcommand{\Pt}[2]{\mathcal{P}_{#1}^{#2}}

\newcommand{\grad}{{\rm rad}_{{\rm P},k}}

\begin{document}
\title{Generalized power domination in WK-Pyramid Networks}
\author{Seethu Varghese \footnote{E-mail: seethu333@gmail.com}\, \, \, \,   A. Vijayakumar\footnote{E-mail: vambat@gmail.com}
\medskip
\\
Department of Mathematics\\ Cochin University of Science and Technology, Cochin-682022, India
}\date{}
 \maketitle
 \begin{abstract}
The notion of power domination arises in the context of monitoring an electric power system with as few phase measurement units as possible. The $k-$power domination number of a graph $G$ is the minimum cardinality of a $k-$power dominating set ($k-$PDS) of $G$. In this paper, we determine the $k-$power domination number of WK-Pyramid networks, $WKP_{(C,L)}$, for all positive values of $k$ except for $k=C-1, C \geq 2$, for which we give an upper bound. The $k-$propagation radius of a graph $G$ is the minimum number of propagation steps needed to monitor the graph $G$ over all minimum \kPDS. We obtain the $k-$propagation radius of $WKP_{(C,L)}$ in some cases. 
 \end {abstract}
\ni \textbf{Mathematics Subject Classification (2010):} 05C69, 94C15\\
\ni \textbf{Keywords:} power domination, electrical network monitoring, domination, WK-Pyramid network
\section {Introduction}

Power domination is a variation of domination introduced in \cite{bamiboad-93} to address the problem of monitoring electrical networks with {\em phasor measurement units}. It was described as a graph parameter in \cite{hahehehe-02}. Let $G =(V(G), E(G))$ be a graph that represents an electric power system, where a vertex $v\in V(G)$ represents an electrical node and an edge $e\in E(G)$ represents a transmission line joining two electrical nodes. For a set $S\subseteq V(G)$, the \textit{closed neighbourhood} of $S$ is the union of the closed neighbourhoods of its elements$\colon N_{G}[S] =\bigcup_{v\in S}N_{G}[v]$ and $< S >$ denotes the subgraph induced by $S\subseteq V(G)$. A vertex $v$ in a graph is said to dominate its closed neighbourhood $N_G[v]$. A subset $S\subseteq V(G)$ of vertices is a dominating set if $N_G[S]=V(G)$, that is if every vertex in the graph is dominated by some vertex of $S$. The minimum size of a dominating set in a graph $G$ is called its \textit{domination number}, denoted by $\gamma(G)$. 

In power domination, there is an additional propagation behaviour. Initially, a set $S$ is said to monitor its closed neighbourhood, like in domination. Then, every vertex that is the only unmonitored neighbour of a monitored vertex gets monitored.  This possibility of propagation conveys the capacity of deducting the status of a node in an electrical network by applying Kirchhoff laws. It gives to power domination a very different flavour since a vertex may then eventually monitor another vertex far apart as in the case of paths. 

The power domination problem was proved to be NP-complete for bipartite graphs and chordal graphs~\cite {hahehehe-02}. Linear-time algorithms for this problem were known for trees~\cite{hahehehe-02}, interval graphs~\cite{lile-05} and block graphs~\cite{xukashzh-06}. 
Upper bounds for the power domination number were studied in~\cite{zhkach-06, ZK-07} and closed formulae for the power domination number were also determined for some graphs~\cite{dohe-06, domoklsp-08, FSV-11}.

Power domination was then generalized in \cite{chdomora-12} by adding the possibility of propagating up to $k$ vertices, $k$ a non-negative integer. Formally, the set of monitored vertices is then described with the following definition from \cite{chdomora-12,dohelomora-13}, inspired by what was proposed in \cite{aaz-10}:
\bd[Monitored vertices]
Let $G$ be a graph, $S \subseteq V(G)$ and $k\ge 0$. The sets $\big(\Pt {G,k} i(S) \big)_{i\ge 0}$ of {\em vertices monitored by $S$ at step $i$} are defined as follows:
\begin{quote}
$\Pt {G,k} 0 (S) = N_{G}[S]$, and \\ 
$\Pt {G,k} {i+1} (S) = \bigcup \{ N_{G}[v]\colon v\in \Pt {G,k} i (S)$
    such that $\big|N_{G}[v]\setminus \Pt {G,k} i (S) \big| \le k\}$.
\end{quote}
\ed

The second part represents the propagation rule. Since $\Pt {G,k} {i} (S)$ is always a union of neighbourhoods, $\Pt {G,k} {i} (S) \subseteq \Pt {G,k} {i+1} (S)$. If $\Pt {G,k} {i_{0}}(S)=\Pt {G,k} {i_{0}+1}(S)$ for some $i_0$, then $\Pt {G,k} {j} (S)=\Pt {G,k} {i_{0}}(S)$ for any $j\geq i_0$. We thus define $\Pt {G,k} {\infty}(S)=\Pt {G,k} {i_{0}}(S)$. When the graph $G$ is clear from the context, we will simplify the notations to $\Pt{k}{i}(S)$ and $\Pt{k}{\infty}(S)$.
 
\bd\cite{chdomora-12}
A $k-${\em power dominating set} of $G$ ($k-$PDS) is a set $S\subseteq V(G)$ such that $\mathcal{P}_{G,k}^{\infty}(S)=V(G)$. The \textit{$k-$power domination number}, $\gamma_{P,k}(G)$, of $G$ is the least cardinality of a $k-$power dominating set of $G$. A $\gpk(G)$-{\em set} is a \kPDS\ in $G$ of cardinality $\gpk(G)$. 
\ed

Generalized power domination reduces to the usual power domination when $k=1$ and to the domination when $k=0$. In~\cite{chdomora-12}, Chang et al. extended several known results for power domination to $k-$power domination. 
They gave sharp upper bounds for the generalized power domination number of connected graphs and of connected claw-free $(k+2)$-regular graphs. In ~\cite{dokl-14}, the authors introduced the $k$-{\em propagation radius} of a graph $G$, motivated from the studies in~\cite{aaz-10}, as a way to measure the efficiency of a minimum \kPDS. It gives the minimum number of propagation steps needed  to monitor the entire graph $G$ over all $\gpk(G)-$sets. They investigated the relationship between propagation radius and the radius of a graph and also computed the propagation radius of Sierpi\'{n}ski graphs.

\bd\cite{dokl-14}
The radius of a \kPDS\ is defined by 
$$\grad(G,S) = 1 + \min \{i:\ \Pt {G,k} {i}(S) = V(G)\}\,.$$
The $k-${\em propagation radius} of the graph as defined in \cite{dokl-14} can be expressed as
$$\grad(G) = \min \{\grad(G,S),\ S\ \text{is a}\ k\text{-PDS\ of}\ G,\ |S| = \gpk(G) \}\,.$$
\ed 

The WK-Pyramid network, an interconnection network based on the WK-recursive mesh~\cite{vesa-88}, was introduced in~\cite{hosa-05} for massively parallel computers. 
It has interesting topological characteristics making it suitable for utilization as the base topology of large scale multi-computer systems. It eliminates some drawbacks of the conventional pyramid network, stemming from the fact that the connections within the layers of this network form a WK-recursive mesh. It is of much less network cost than the hypercube, $k-$ary $n-$cube and WK-Recursive networks. It also has small average distance and diameter, large connectivity and high degree of scalability and expandability. Because of the desirable properties of this network, it is suitable for medium or large sized networks and also a best alternative for mesh and traditional pyramid interconnection topologies.

For $C, L\in \mathbb N$, let $[C]=\{1,\hdots, C\}$, $[C]_{0}=\{0,\hdots, C-1\}$ and $[C]_0^{L-2}=\{a_{L-2}a_{L-3}\hdots a_{1}\colon a_{i}\in [C]_{0}$ for all $i\}$.
An $L-$level {\em WK-Recursive mesh}~\cite{hosa-05}, denoted by $WK_{(C,L)}$, consists of a set of vertices \\$V(WK_{(C,L)})=\{(a_{L}a_{L-1}\hdots a_{1})\colon a_{i}\in [C]_0$ for $i\in [L]\}$. The vertex with address $(a_{L}a_{L-1}\hdots a_{1})$ is adjacent 
\begin{enumerate}
\item to all the vertices with addresses $(a_{L}a_{L-1}\hdots a_{2}a_{j})$ such that $a_{j}\in [C]_0$, $a_{j}\neq a_{1}$ and
\item to a vertex $(a_{L}a_{L-1}\hdots a_{j+1}a_{j-1}(a_{j})^{j-1}),$ if there exists one $j$ such that $2\leq j\leq L, a_{j-1}= a_{j-2}=\hdots=a_{1}$ and $a_{j}\neq a_{j-1}$. 
\end{enumerate}
The notation $(a_{j})^{j-1}$ denotes that the term $a_{j}$ is repeated $j-1$ times. Vertices of the form $(\stackrel{L \text{ times}}{\overbrace{a\hdots a}})$ are called extreme vertices of $WK_{(C,L)}$. Clearly, $WK_{(C,L)}$ contains $C$ extreme vertices of degree $C-1$ and all the other vertices are of degree $C$. Note that $WK_{(C,1)}\cong K_{C}$ $(C\geq1)$, $WK_{(1,L)}\cong K_{1}$ and $WK_{(2,L)}\cong P_{2^L}$ $(L\geq1)$, where $K_{C}$ and $P_{L+1}$ denote the complete graph of order $C$ and the path of order $2^L$, respectively.

A {\em WK-Pyramid network}~\cite{hosa-05}, denoted by $WKP_{(C,L)}$, consists of a set of vertices\\$V(WKP_{(C,L)})=\{(r,(a_{r}a_{r-1}\hdots a_{1}))\colon r\in [L], a_{i}\in [C]_0$ for $i\in [r]\}\cup \{(0,(1))\}$. A vertex with addressing scheme $(r,(a_{r}a_{r-1}\hdots a_{1}))\in V(WKP_{(C,L)})$ is called a vertex at level $r$. The part $(a_{r}a_{r-1}\hdots a_{1})$ of the address determines the address of a vertex within the WK-recursive mesh at level $r$. The vertex $(0,(1))$ is adjacent to every vertex in level 1. 
A vertex with address $(r,(a_{r}a_{r-1}\hdots a_{1}))$ at level $r > 0$ is adjacent
\begin{enumerate}
\item to vertices $(r,(a_{r}a_{r-1}\hdots a_{2}a_{j}))\in V(WKP_{(C,L)}),$ for $a_{j}\in [C]_0$, $a_{j}\neq a_{1}$  
\item to a vertex with address schema $(r,(a_{r}a_{r-1}\hdots a_{j+1}a_{j-1}(a_{j})^{j-1}))$, if there exists one $j$ such that $2\leq j\leq L, a_{j-1}= a_{j-2}=\hdots=a_{1}$ and $a_{j}\neq a_{j-1}$
\item to vertices $(r+1,(a_{r}a_{r-1}\hdots a_{2}a_{1}a_{j})),$ for $a_{j}\in [C]_0$, in level $r+1$ and
\item to a vertex $(r-1,(a_{r}a_{r-1}\hdots a_{2}))$, in level $r-1$.  
\end{enumerate} 

\begin{figure}[h]\center
  \includegraphics[width=8cm]{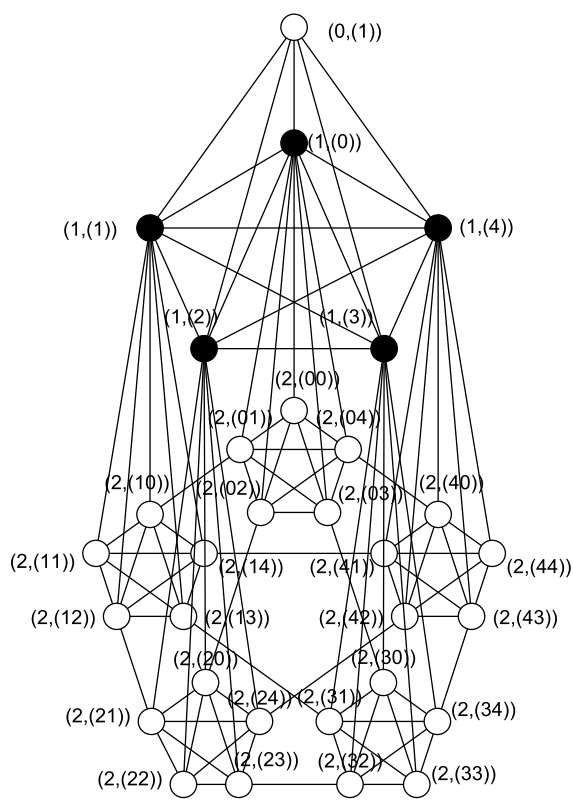}\\[-0.3cm]
  \caption{The graph $WKP_{(5,2)}$ }\label{pic1}
\end{figure}

All the vertices in level $r>0$ of $WKP_{(C,L)}$ induce a WK-recursive mesh $WK_{(C,r)}$. Note that $WKP_{(C,1)}\cong K_{C+1}$ $(C\geq1)$, $WKP_{(1,L)}\cong P_{L+1}$ $(L\geq1)$. Fig.~\ref{pic1} shows the graph $WKP_{(5,2)}$. Vertices of the form $(r,(\stackrel{r \text{ times}}{\overbrace{a\hdots a}}))$ are called the extreme vertices of $WKP_{(C,L)}$. The vertex $(0,(1))$ has degree $C$ and at any level except the $L^{th}$ level, the extreme vertices  are of degree $2C$ and the other vertices are of degree $2C+1.$ In the $L^{th}$ level, the extreme vertices have degree $C$ and the other vertices have degree $C+1$.

\begin{figure}[h]\center
  \includegraphics[width=4cm]{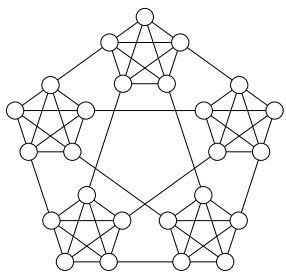}\\[-0.3cm]
  \caption{The induced subgraph $G_{w}^{5,L}$ of $WKP_{(5,L)}$}\label{pic2}
\end{figure}

For $w\in [C]_0^{L-2}$, let $V_w^{C,L}=\{(L,(wij))\in WKP_{(C,L)}\colon i,j\in [C]_{0}\}$ and $G_{w}^{C,L} = \\< V_w^{C,L} >$, i.e. $G_{w}^{C,L}$ is the induced subgraph in level $L$ of $WKP_{(C,L)}$. In fact, $G_{w}^{C,L}$ is isomorphic to $WK_{(C,2)}$ for any $w\in [C]_0^{L-2}$ and any $L\geq3$ (Fig.~\ref{pic2}).%


The $k-$power domination number is known only for a few nontrivial families of graphs. In this paper, we determine the $k-$power domination number of WK-Pyramid network for all positive values of $k$ except for $k=C-1, C \geq 2$, for which we give an upper bound. This is the first network class with the pyramid structure for which the $k-$power domination number is studied. We also obtain the $k-$propagation radius of $WKP_{(C,L)}$ in some cases.

\section{$\boldsymbol{k-}$power domination number of WK-Pyramid network}
\bt\label{th:k>C}
Let $C, L, k\geq 1$. If $C=1$ or $L=1$ or $k\geq C$, then $\gpk(WKP_{(C,L)})=1$ .
\et
\bp
Recall that $WKP_{(C,1)}\cong K_{C+1}$ $(C\geq1)$ and that $WKP_{(1,L)}\cong P_{L+1}$ $(L\geq1)$. Hence $\gpk(G)=1$ for these graphs $G$. If $k\geq C$, then take $S=\{(0,(1))\}$. It monitors the vertices in level 1. Since each vertex in level $r$ has exactly $C$ neighbours in its successive level $r+1$, once the level $r$ is monitored, the vertices in level $r+1$ get monitored by propagation. This propagation goes on till level $L$ and hence $S$ is a $k-$PDS of $WKP_{(C,L)}$.   
\ep
We now consider the case when $k\le C-1$. We begin with the computation of $\gpk$ for $L=2$. Let $V^1$ and $V^2$ denote the set of vertices of $WKP_{(C,2)}$ in levels 1 and 2 respectively. Let $Q_i$ denote a $C-$clique induced by the set of vertices $\{(2,(ij))\colon j\in [C]_0\}$ for some $i$. We first obtain the following upper bound.
\bl \label{l:1}
For $C\geq3$ and $k\in [C-1]$, $\gpk(WKP_{(C,2)})\leq C-k.$
\el
\bp
Let $S=\{(1,(i))\colon k\leq i\leq C-1\}$. (The vertices in $S$ are coloured black in Fig.~\ref{pic1}.)\\
Then $\mathcal{P}_k^{0}(S)=\begin {cases}(1,(j))&;j\in [C]_0\\
                          (2,(ij))&; k\leq i\leq C-1, j\in [C]_0\\
                          (0,(1)) \end{cases}$\\ 
$\mathcal{P}_k^{1}(S)=\mathcal{P}_k^{0}(S)\cup\{(2,(ij))\colon i\in [k]_0, k\leq j\leq C-1\}$ and\\ 
$\mathcal{P}_k^{2}(S)=\mathcal{P}^{1}(S)\cup\{(2,(ij))\colon i, j \in [k]_0\} =V(WKP_{(C,2)})$.\\
Hence $S$ is a $k-$PDS, which implies $\gpk(WKP_{(C,2)})\leq \left|S\right|= C-k.$
\ep
\bl \label{l:2}
For $C\geq3$ and $k\in [C-2]$, $\gpk(WKP_{(C,2)})\geq C-k.$
\el
\bp
Let $S$ be a minimum \kPDS\ of $WKP_{(C,2)}$. We may assume that $S\subseteq V^1\cup V^2$.

\textbf{Claim:} $\left|S\cap(V^{1}\cup V^{2})\right|\geq C-k.$

Suppose on the contrary that $\left|S\cap(V^{1}\cup V^{2})\right|\leq C-k-1.$ We consider the case when $S$ contains vertices from both $V^{1}$ and $V^{2}.$ 
Assume first that $\left|S\cap(V^{1}\cup V^{2})\right|= C-k-1$ and that $S$ contains a vertex $(1,(i'))\in V^{1}$ and the remaining $C-k-2$ vertices from the $C-$cliques $Q_{i_1},\hdots, Q_{i_{C-k-2}}$, where $i'\neq i_\ell, \ell\in [C-k-2]$ such that each of these $C-$cliques contains exactly one vertex in $S$. Let $Q_{\ell}$ be an arbitrary clique that does not contain any vertex of $S$, where $\ell\neq i'$. Let $X=\{(2,(\ell i'))\}\cup\{(2,(\ell i_1)),\hdots,(2,(\ell i_{C-k-2}))\}$. Then $\mathcal{P}_k^{1}(S)\cap Q_{\ell}= X$. This holds for every $l\in L'=[C]_{0}\setminus \{i', i_1,\hdots,i_{C-k-2}\}$. Thus the set of vertices $J=\{(2,(\ell\ell'))\colon \ell\in L', \ell'\in L'\}$ has an empty intersection with $\mathcal{P}_k^{1}(S)$. Since every vertex in $WKP_{(C,L)}\setminus J$ has either 0 or $k+1$ neighbours in $J$, no vertex from this set $J$ may get monitored later on, which is a contradiction. Assume next that $\left|S\cap(V^{1}\cup V^{2})\right|< C-k-1$ or that $S$ intersects some $C-$clique $Q_{i}$ in more than one vertex. Then we can analogously conclude that not all vertices of $Q_{\ell}$ will be monitored. Now, the case when $S\cap V^{1} = \phi$ or $S\cap V^{2} = \phi$ can be proved in a similar manner. Hence the claim.

Therefore $\gpk(WKP_{(C,2)})=|S|=|S\cap(V^{1}\cup V^{2})|\geq C-k.$
\ep

\bt \label{th:wkpc2}
For $C\geq2$ and $k\in [C-1]$, $\gpk(WKP_{(C,2)})= C-k.$
\et
\bp
Clearly, $\gamma_{\rm P,1}(WKP_{(2,2)})= 1.$ Let $C\geq3$. For $k= C-1$, any vertex in level 1 forms a $k-$PDS of $WKP_{(C,2)}$. For $k\in [C-2]$, the result follows from Lemmas~\ref{l:1} and~\ref{l:2}.  
\ep

\bl \label{l:3}
For $C\geq3, L\geq3$ and $k\in [C-2]$, $\gpk(WKP_{(C,L)})\leq (C-k-1)C^{L-2}.$ 
\el

\bp

In $WKP_{(C,L)}$, the vertices in the $L^{th}$ level induce $WK_{(C,L)}$ which is hamiltonian~\cite{hosa-05}. Also, by contracting each of the subgraphs $G_{w}^{C,L}$ into a single vertex, the graph induced by the vertices in level $L$ is isomorphic to $WK_{(C,L-2)}$. Hence, in level $L$ of $WKP_{(C,L)}$, we can arrange the subgraphs of the form $G_{w}^{C,L}$ into a cycle such that there exists exactly one edge between the consecutive subgraphs. We now construct a set $S$ in such a way that corresponding to each subgraph $G_{w}^{C,L}$ in level $L$, $S$ contains one vertex from the neighbour set of $G_{w}^{C,L}$ in level $L-1$ (which induces a clique) and $C-k-2$ additional vertices from $G_{w}^{C,L}$.

Let $w',w''\in [C]_0^{L-2}$. 
Let $G_{w}^{C,L},G_{w'}^{C,L}$ and $G_{w''}^{C,L}$ be consecutive subgraphs in the selected hamiltonian order. Let $xx'$ be the edge between $G_{w}^{C,L}$ and $G_{w'}^{C,L}$, where $x\in G_{w}^{C,L}, x'\in G_{w'}^{C,L}$ and let $y'y''$ be the edge between $G_{w'}^{C,L}$ and $G_{w''}^{C,L}$, where $y'\in G_{w'}^{C,L}, y''\in G_{w''}^{C,L}$. Denote $x=(L,(wii))$ and $y'=(L,(w'jj))$ for some $i$ and $j$, $i\neq j$. Let $S$ contain the vertex $(L-1,(w'j))$ (which is the neighbour of $y'$ in the $(L-1)^{th}$ level) and $C-k-2$ additional vertices from $G_{w'}^{C,L}$ such that no two lying in the same $C-$clique in $G_{w'}^{C,L}$ and no one lying in the $C-$clique $H$ in $G_{w'}^{C,L}$ that contains $x'$. Also $S\cap Q=\phi$, where $Q$ is the $C-$clique in $G_{w'}^{C,L}$ that contains $y'$. Now, do this in parallel for all the corresponding subgraphs. In particular, the vertex $(L-1,(wi))$ in the $(L-1)^{th}$ level corresponding to the vertex $x$ is put into $S$ when considering $G_{w}^{C,L}$. Thus $C-k$ vertices of $H$ lie in $\mathcal{P}_k^{1}(S)$: one of these vertices is $x'$, the other $C-k-1$ are those vertices of $H$ that have a neighbour in the $C-$cliques in $G_{w'}^{C,L}$ that contain $C-k-2$ vertices of $S$ and that have a neighbour in the $C-$clique $Q$ in $G_{w'}^{C,L}$. Also the neighbour of $H$ in the $(L-1)^{th}$ level belongs to $\mathcal{P}_k^{0}(S)$ since $(L-1,(wi))\in S$. Hence the remaining $k$ vertices of $H$ lie in $\mathcal{P}_k^{2}(S)$ and it is straightforward to check that all the vertices of $G_{w'}^{C,L}$ lie in $\mathcal{P}_{k}^{\infty}(S)$. In a similar way, every vertex in the $L^{th}$ level is monitored. We know that, for any $w$, the neighbours of $G_{w}^{C,L}$ in the $(L-1)^{th}$ level induce a $C-$clique. By the construction of $S$, each $C-$clique in the $(L-1)^{th}$ level contains a vertex in $S$. Thus we get that all the vertices in levels $L-1$ and $L-2$ belong to $\mathcal{P}_k^{0}(S)$. Now, since each vertex in level $L-2$ has exactly one neighbour in its preceeding level, vertices in the $(L-3)^{rd}$ level are monitored by propagation. This propagation continues to the preceeding levels and hence the whole graph gets monitored. Thus we conclude that $S$ is a $k-$PDS. Since each subgraph $G_{w}^{C,L}$ contains $C-k-1$ vertices of $S$, $\left|S\right|\leq (C-k-1)C^{L-2}.$     
\ep

\bl \label{l:4}
For $C\geq3, L\geq3$ and $k\in [C-2]$, $\gpk(WKP_{(C,L)})\geq (C-k-1)C^{L-2}.$ 
\el
\bp
Let $S$ be a minimum \kPDS\ of $WKP_{(C,L)}$ and $w\in [C]_0^{L-2}.$  
Denote $V_{w}^{C,L-1}=\{(L-1,(wi))\in WKP_{(C,L)}\colon i\in [C]_0\}$. 

\textbf{Claim:} $\left|S\cap(V_w^{C,L}\cup V_w^{C,L-1})\right|\geq C-k-1.$

Suppose on the contrary that $\left|S\cap(V_w^{C,L}\cup V_w^{C,L-1})\right|\leq C-k-2.$ Consider the case when $S\cap V_w^{C,L-1}=\phi.$ Then $|S\cap V_w^{C,L}|\leq C-k-2.$ Assume first that $\left|S\cap V_w^{C,L}\right|=C-k-2$. Let $H_{i}$ be a $C-$clique in $G_{w}^{C,L}$, i.e. $H_{i}$ is induced by the set of vertices $\{(L,(wij))\in WKP_{(C,L)}\colon j\in [C]_0\}$ for some $i$. Assume that $S$ has exactly one vertex in  $C-$cliques $H_{i}$ for $i\in \{i_1,\hdots,i_{C-k-2}\}$. Then $S\cap V(H_{i'})=\phi$ holds for other $k+2$ coordinates $i'$. Let $H_{\ell}$ be an arbitrary such clique in $G_{w}^{C,L}$ that does not contain any vertex of $S$. Let $X=\{(L,(w\ell i_1)),\ldots,(L,(w\ell i_{C-k-2}))\}\cup \{(L,(w\ell \ell))\}$. Then $\mathcal{P}_k^{1}(S)\cap H_{\ell}\subseteq X$. This holds for every $\ell\in L'=[C]_{0} \setminus\{i_1,\hdots,i_{C-k-2}\}$. Thus the set of vertices $\{(L,(w\ell\ell'))\colon \ell\in L', \ell'\in L', \ell\neq \ell'\}$ has an empty intersection with $\mathcal{P}_k^{1}(S)$. Since every vertex in $WKP_{(C,L)}$ has either 0 or $k+1$ neighbours in this set, no vertex from this set may get monitored later on, a contradiction. Assume next that $\left|S\cap V_w^{C,L}\right|< C-k-2$ or that $S$ intersects some $C-$clique $H_{i}$ in more than one vertex. Then we can analogously conclude that not all vertices of $H_{\ell}$ will be monitored. Thus the case when $S\cap V_w^{C,L-1}=\phi$ is not possible. 

Now suppose that $S \cap V_w^{C,L-1} \neq \phi$. Assume first that $\left|S\cap(V_w^{C,L} \cup V_w^{C,L-1})\right|= C-k-2$ and that $S$ contains a vertex $(L-1,(wi'))\in V_w^{C,L-1}$ and the remaining $C-k-3$ vertices from the $C-$cliques $H_{i_1},\hdots, H_{i_{C-k-3}}$, where $i'\neq i_\ell, \ell\in [C-k-3]$ such that each of these $C-$cliques contains exactly one vertex in $S$. Let $H_{\ell}$ be an arbitrary clique in $G_{w}^{C,L}$ that does not contain any vertex of $S$, where $\ell\neq i'$. Let $X=\{(L,(w\ell i'))\}\cup \{(L,(w\ell\ell))\}\cup\{(L,(w\ell i_1)),\hdots,(L,(w\ell i_{C-k-3}))\}$. Then $\mathcal{P}_k^{1}(S)\cap H_{\ell}\subseteq X$. This holds for every $l\in L'=[C]_{0} \setminus \{i', i_1,\hdots,i_{C-k-3}\}$. Thus the set of vertices $\{(L,(w\ell\ell'))\colon \ell\in L', \ell'\in L', \ell\neq \ell'\}$ has an empty intersection with $\mathcal{P}_k^{1}(S)$. Since every vertex in $WKP_{(C,L)}$ has either 0 or $k+1$ neighbours in this set, no vertex from this set may get monitored later on, which is a contradiction. Assume next that $\left|S\cap(V_w^{C,L}\cup V_w^{C,L-1})\right|< C-k-2$ or that $S$ intersects some $C-$clique $H_{i}$ in more than one vertex. Then we can analogously conclude that not all vertices of $H_{\ell}$ will be monitored. Hence the claim. Therefore, $\left|S\cap(V_w^{C,L}\cup V_w^{C,L-1})\right|\geq C-k-1$, i.e. $\left|S\cap (V(G_{w}^{C,L})\cup N_{L-1}(G_{w}^{C,L}))\right|\geq C-k-1$, where $N_{L-1}(G_{w}^{C,L})$ is the set of neighbours of $G_{w}^{C,L}$ in the $(L-1)^{th}$ level. Hence corresponding to each $G_{w}^{C,L}$  in the $L^{th}$ level, we get at least $C-k-1$ vertices in $S$. 

Hence $\left|S\right|
\geq \sum_{w\in[C]_0^{L-2}}(C-k-1)=(C-k-1)C^{L-2}$.        
 \ep
 \bt
For $C\geq3, L\geq3$ and $k\in [C-2]$,\\ $\gpk(WKP_{(C,L)})=(C-k-1)C^{L-2}.$
\et
\bp
Follows from Lemmas~\ref{l:3} and~\ref{l:4}.
\ep 

Thus we have the following consolidated result:\\\\ 
Let $C, L, k\geq 1$. Then 

$\gpk(WKP_{(C,L)})=\begin {cases}1; &C=1$ or $L=1$ or $k\geq C,\\
                       C-k; &L=2, C\geq 2, k\in [C-1],\\
                       (C-k-1)C^{L-2};&L\geq3, C\geq3, k\in [C-2].\end{cases}$\\\\
For $k=C-1$, $C\geq2$ and $L\geq3$ we prove the following upper bound.
\bt
For $C\geq2,$ $L\geq3$, $\gamma_{\rm P,C-1}(WKP_{(C,L)})\leq\left\lceil\frac{L+1}{3}\right\rceil$.
\et
\bp
We consider three cases.
\bca
 $L=3m$
 \eca
 $S=\{\ \bigcup^{m}_{i=1} (3i-1,(0)^{3i-1})\}\cup \{(0,(1))\}$.
 
 Here, $\left|S\right|=m+1$. Also, $\left\lceil\frac{L+1}{3}\right\rceil=\left\lceil\frac{(3m)+1}{3}\right\rceil=m+1$.
  \bca
 $L=3m+1$
 \eca
 $S=\{\ \bigcup^{m}_{i=1} (3i,(0)^{3i})\}\cup \{(1,(0))\}$.
 
 Here, $\left|S\right|=m+1$. Also, $\left\lceil\frac{L+1}{3}\right\rceil=\left\lceil\frac{(3m+1)+1}{3}\right\rceil=m+1$.
  \bca
 $L=3m+2$
 \eca
 $S=\{\ \bigcup^{m+1}_{i=1} (3i-2,(0)^{3i-2})\}.$
 
 Here, $\left|S\right|=m+1$. Also, $\left\lceil\frac{L+1}{3}\right\rceil=\left\lceil\frac{(3m+2)+1}{3}\right\rceil=m+1$.\\\\ 
 In each case, $\mathcal{P}_{C-1}^{\infty}(S)=V(WKP_{(C,L)})$ and thus $S$ is a $k-$PDS of order $\left\lceil \frac{L+1}{3}\right\rceil$. Hence $\gamma_{\rm P,C-1}(WKP_{(C,L)})\leq\left\lceil \frac{L+1}{3}\right\rceil$. 
\ep
\section{Propagation radius of WK-Pyramid network}
In this section, we determine the $k-$propagation radius of $WKP_{(C,L)}$ for $C\ge 1$ and $L=1,2$.  If $L=1$, the graph is a complete graph and the propagation radius is 1. If $C=1, \grad(WKP_{(1,L)})=\grad(P_{L+1})= \left\lfloor \frac{L+1}{2}\right\rfloor$.

\bl\label{lem:rad}
Let $C\ge 3$ and $k \in [C-1]$ and $S$ be a minimum \kPDS\ of $WKP_{(C,2)}$. Then $S\cap V^1\neq \phi.$
\el
\bp
Suppose that $S\cap V^{1} = \phi$. Consider the case when $(0,(1))\notin S$. Then by Theorem~\ref{th:wkpc2}, $|S\cap V^{2}|=C-k$. Assume first that $S$ has exactly one vertex in $C-$cliques $Q_{i}$ for $i\in \{{i_{1}},\hdots,{i_{C-k}}\}$. Then $S\cap V(Q_{i'})=\phi$ for $k$ coordinates $i'$. Let $Q_{\ell}$ be an arbitrary such subgraph. Let $X=\{(2,(\ell i_1)),\hdots,(2,(\ell i_{C-k}))\}$. Then $\mathcal{P}_k^{1}(S)\cap V(Q_{\ell})= X$ and $\mathcal{P}_k^{1}(S)\cap <V^{1}>= \{(1,{i_{1}}),\hdots,(1,{i_{C-k}})\}$. This holds for any $\ell\in L=[C]_{0}\setminus \{{i_{1}},\hdots,{i_{C-k}}\}$. Therefore the set of vertices $K=\{(2,(ij))\colon i,j\in L\}\cup \{(1,(i))\colon i\in L\}\cup \{(0,(1))\}$ has an empty intersection with $\mathcal{P}_k^{1}(S)$. Since every vertex of $WKP_{(C,2)}\setminus K$ has either 0 or $k+1$ neighbours in $K$, no vertex from this set may get monitored later on, a contradiction. The case when $(0,(1))\in S$ or that $S$ intersects some $Q_{i}$ in more than one vertex can be proved analogously.             
\ep
\bt
Let $C\ge 2$ and $k\ge 1$. Then\\
$\grad(WKP_{(C,2)})=\begin {cases}2; &k\geq C,\\
                       3; &k\in [C-1].\end{cases}$\\
\et
\bp
For $k\geq C$, $\gpk(WKP_{(C,2)})=1$, by Theorem~\ref{th:k>C} and observe that $\gamma(WKP_{(C,2)})>1$. Therefore, $\grad(WKP_{(C,2)})\ge 2$ (see~\cite[Proposition 4.1]{dokl-14}). And, for the set $S=\{(0,(1))\}$, we get that $\mathcal{P}_k^{0}(S) = S \cup V^{1}$ and $\mathcal{P}_k^{1}(S) = V(WKP_{(C,2)})$. Now let $k\in [C-1]$. For $C=2$, the result easily follows. Let $C\ge 3$. By Theorem~\ref{th:wkpc2}, $\gpk(WKP_{(C,2)})=C-k$ and therefore by Lemma~\ref{lem:rad}, $|S\cap V^{2}|\le C-k-1$ for every minimum \kPDS\ $S$. Then there exists at least $k+1$ $C-$cliques $Q_{i}$ not containing any vertex of $S$. Let $Q_{i'}$ be an arbitrary clique such that $S\cap V(Q_{i'})=\phi$ and $(1,(i'))\notin S$. We prove that the vertex $(2,(i'i'))$ is not in $\mathcal{P}_k^{1}(S)$. Clearly, $(2,(i'i'))\notin \mathcal{P}_k^{0}(S)$. Moreover, $|V(Q_{i'})\cap \mathcal{P}_k^{0}(S)|\le C-k-1$ and $|V(Q_{i'})\setminus \mathcal{P}_k^{0}(S)|\ge k+1$. Therefore any neighbour of $(2,(i'i'))$ in $Q_{i'}$ is adjacent to more than $k$ unmonitored vertices preventing any propagation to this vertex at this step. Also, since $(1,(i'))$ has more than $k$ unmonitored vertices as its neighbours, $(2,(i'i'))$ cannot be monitored by $(1,(i'))$ at this stage. Hence $\grad(WKP_{(C,2)})\ge3$. Also, by Lemma~\ref{l:1}, $\grad(WKP_{(C,2)})\le3$. 
\ep

\textbf{Note:} For $C,L\ge 3$ and $k\ge 1$, by observing the propagation behaviour described in the proof of Theorem~\ref{th:k>C} and Lemma~\ref{l:3}, one can obtain that $\grad(WKP_{(C,L)})\le L$ if $k\ge C$ and $\grad(WKP_{(C,L)})\le\max\{5, L-1\}$ if $k\in [C-2]$.

\section{Conclusion}
In this paper, we have determined the $k-$power domination number of WK-Pyramid networks, $WKP_{(C,L)}$, for all positive values of $k$ except for $k=C-1, C \geq 2$, for which we give an upper bound. We also obtain the $k-$propagation radius of $WKP_{(C,L)}$ in some cases. The $k-$power domination number of other pyramid networks such as grid pyramids, torus pyramids can be studied in future.\\\\

\textbf{Acknowledgement:} The first author is supported by Maulana Azad National Fellowship (F1- 17.1/2012-13/MANF-2012-13-CHR-KER-15793) of the University Grants Commission, India.

\end{document}